\documentclass{amsart}
\usepackage{amsfonts}
\usepackage{graphicx}
\usepackage{amsmath}
\usepackage{geometry}

\newtheorem{theorem}{Theorem}

\newtheorem{definition}[theorem]{Definition}
\newtheorem{example}[theorem]{Example}

\newtheorem{lemma}[theorem]{Lemma}
\newtheorem{notation}[theorem]{Notation}

\newtheorem{remark}[theorem]{Remark}

\input{tcilatex}
\begin{document}
\title{Factorization in Quantum Planes}
\author{Romain Coulibaly and Kenneth Price}
\address{Department of Mathematics\\
University of Wisconsin Oshkosh\\
800 Algoma Boulevard\\
Oshkosh, WI 54901}
\email{angecoul@yahoo.com or pricek@uwosh.edu}
\date{}
\maketitle

\begin{abstract}
These results stem from a course on ring theory. \ Quantum planes are rings
in two variables $x$ and $y$ such that $yx=qxy$ where $q$ is a nonzero
constant. \ When $q=1$ a quantum plane is simply a commutative polynomial
ring in two variables. \ Otherwise a quantum plane is a noncommutative ring.

Our main interest is in quadratic forms belonging to a quantum plane. \ We
provide necessary and sufficient conditions for quadratic forms to be
irreducible. \ We find prime quadratic forms and consider more general
polynomials. \ Every prime polynomial is irreducible and either central or a
scalar multiple of $x$ or of $y$. \ Thus there can only be primes of degree
2 or more when $q$ is a root of unity.
\end{abstract}

\section{Introduction}

Throughout $F$ denotes a commutative field of characteristic zero. \ We
often let $F$ equal $\mathbb{Q}$, the rational number field, or $\mathbb{R}$%
, the real number field. \ Let $\mathbb{N}$ denote the set of nonnegative
integers and let $\mathbb{N}^{2}$ denote the set of ordered pairs with
entries from $\mathbb{N}$.

This work was motivated by an abstract algebra course taught by Kenneth
Price. \ Many of the results were subsequently obtained by Romain Coulibaly
while still an undergraduate student. \ Our interest stems from the course
topics on principal ideal domains (PID's) and unique factorization domains
(UFD's). \ For example, $F\left[ x\right] $ is a PID. \ Every PID\ is a UFD
but $F\left[ x,y\right] $ is a UFD which is not a PID. \ A standard result,
which generalizes Euclid's Lemma, is that irreducible elements of a UFD are
prime (see \cite[Lemma 9.2.1]{BB}).

Background material is covered in section \ref{Pq def section}. \ Next we
begin the study of quadratic forms and introduce a quantum analog of the
discriminant in section \ref{Irreducible Sec}. \ We give conditions for a
quadratic form to be reducible. \ We provide examples of reducible quadratic
forms which have two distinct factorizations into irreducible polynomials
when $F=\mathbb{R}$.

Prime elements of general noncommutative rings are discussed in section \ref%
{Noncommutative Rings Section}. \ Prime elements of quantum planes are
covered in section \ref{Prime Poly Sec}. \ Examples of prime quadratic forms
are provided. \ We show any prime polynomial must be central. \ A final
example is an irreducible and central polynomial which is not prime.

\section{The Quantum Plane \label{Pq def section}}

\begin{notation}
Let $q\in F\backslash \left\{ 0\right\} $ be arbitrary.\ \ The \emph{quantum
plane}, denoted $\mathcal{O}_{q}\left( F^{2}\right) $, is a ring generated
by variables $x,y$ whose product is governed by equation \ref{commutation
equation}. 
\begin{equation}
yx=qxy  \label{commutation equation}
\end{equation}
\end{notation}

Quantum planes are well-known to noncommutative ring theorists and\ were
first introduced as the noncommutative analog of a coordinate ring of a two
dimensional space by Yuri Manin (see \cite[section 4.2.8]{Manin}). \ Our
notation follows \cite[Example 1.1.4]{GB} and \cite[Proposition 1.1.13]{GW},
although we have interchanged $x$ and $y$.

\begin{definition}
Elements of $\mathcal{O}_{q}\left( F^{2}\right) $ are called \emph{%
polynomials}. \ The \emph{terms} of $\mathcal{O}_{q}\left( F^{2}\right) $
have the form $x^{i}y^{j}$ where $i,j$ are nonnegative integers and $%
1=x^{0}y^{0}$.
\end{definition}

We may denote a polynomial by $f\left( x,y\right) $ or simply by $f$. \ The
set of terms is a basis of $\mathcal{O}_{q}\left( F^{2}\right) $. \ For each 
$f\in \mathcal{O}_{q}\left( F^{2}\right) $ there is a unique expression in
equation \ref{polynomial sum} with $\lambda _{\left( i_{1},i_{2}\right) }\in
F$ and only finitely many $\lambda _{\left( i_{1},i_{2}\right) }\neq 0$\ for 
$\left( i_{1},i_{2}\right) \in \mathbb{N}^{2}$.%
\begin{equation}
f=\sum_{\left( i_{1},i_{2}\right) \in \mathbb{N}^{2}}\lambda _{\left(
i_{1},i_{2}\right) }x^{i_{1}}y^{i_{2}}  \label{polynomial sum}
\end{equation}

We often do not specify a value for $q$, but treat it like a variable which
commutes with $x$, $y$, and elements of $F$. \ To simplify the product of
elements of $\mathcal{O}_{q}\left( F^{2}\right) $ we put $x$ to the left of $%
y$ in each term. \ This can be achieved since equation \ref{commutation
equation} extends to $y^{j}x^{i}=q^{j\cdot i}x^{i}y^{j}$ for all $i,j\in 
\mathbb{N}$.

\begin{definition}
Let $f$ be written as in equation \ref{polynomial sum}.

\begin{enumerate}
\item We call $\lambda _{\left( i_{1},i_{2}\right) }$ a \emph{coefficient}
of $f$ for $\left( i_{1},i_{2}\right) \in \mathbb{N}^{2}$.

\item The \emph{degree} of $f$, denoted $\deg \left( f\right) $, is the
largest value of $i_{1}+i_{2}$ for all $\left( i_{1},i_{2}\right) \in 
\mathbb{N}^{2}$ such that $\lambda _{\left( i_{1},i_{2}\right) }\neq 0$.

\item We say $f$ is \emph{homogeneous} if $i_{1}+i_{2}=\deg \left( f\right) $
for all $\left( i_{1},i_{2}\right) \in \mathbb{N}^{2}$ with $\lambda
_{\left( i_{1},i_{2}\right) }\neq 0$.

\item The \emph{support} of $f$ is $\left\{ x^{i_{1}}y^{i_{2}}:\left(
i_{1},i_{2}\right) \in \mathbb{N}^{2}\text{ and }\lambda _{\left(
i_{1},i_{2}\right) }\neq 0\right\} $.
\end{enumerate}
\end{definition}

The proof of Lemma \ref{homogeneous Lemma} is similar to the proof of
Gauss's Lemma provided in most introductory algebra texts, such as \cite[%
Lemma 4.3.3 and Theorem 4.3.4]{BB}.

\begin{lemma}
\label{homogeneous Lemma}Suppose $f,g\in \mathcal{O}_{q}\left( F^{2}\right) $
are such that $fg$ is homogeneous. \ If $q\neq 0$ then $f$ and $g$ are
homogeneous.
\end{lemma}

\begin{proof}
We proceed by contradiction and assume one of $f$ and $g$ is not homogeneous
and $fg$ is homogeneous. Write $f$ as in equation \ref{polynomial sum} and
write $g$ as in equation \ref{gseq} with $\mu _{\left(j_{1},j_{2}\right)
}\in F$ and only finitely many $\mu _{\left( j_{1},j_{2}\right) }\neq 0$\
for $\left( j_{1},j_{2}\right) \in \mathbb{N}^{2}$. 
\begin{equation}
g=\sum_{\left( j_{1},j_{2}\right) \in \mathbb{N}^{2}}\mu _{\left(
j_{1},j_{2}\right) }x^{j_{1}}y^{j_{2}}  \label{gseq}
\end{equation}

Let $\prec $ be the linear ordering on $\mathbb{N}^{2}$ given by $\left(
a_{1},a_{2}\right) \prec \left( b_{1},b_{2}\right) $ if $\left( a_{1}+a_{2}
\right) < \left( b_{1}+b_{2}\right) $ or $\left( a_{1}+a_{2} \right) =
\left( b_{1}+b_{2}\right)$ and $a_{1}<b_{1}$ for $\left( a_{1},a_{2}\right)
,\left( b_{1},b_{2}\right) \in \mathbb{N}^{2}$. \ Choose $\left(
i_{1},i_{2}\right) ,\left( j_{1},j_{2}\right) \in \mathbb{N}^{2}$ as small
as possible with respect to $\prec $\ so that $\lambda _{\left(
i_{1},i_{2}\right) }$ $\neq 0$ and $\mu _{\left( j_{1},j_{2}\right) }\neq 0$%
. \ Then $i_{1}+i_{2}+j_{1}+j_{2} < \deg \left( f\right) +\deg \left(
g\right) $ since one of $f$ and $g$ is not homogeneous.

The coefficient of $x^{i_{1}+j_{1}}y^{i_{2}+j_{2}}$ in $fg$ must be zero
since $fg$ is homogeneous of degree $\deg \left( fg \right) = \deg \left(
f\right) +\deg \left( g\right) $. This yields equation \ref{prod eq} with
the sum over all $\left( i_{1}^{\prime },i_{2}^{\prime }\right) ,\left(
j_{1}^{\prime },j_{2}^{\prime }\right) \in \mathbb{N}^{2}$ such that $%
i_{1}^{\prime }+j_{1}^{\prime }=i_{1}+j_{1}$ and $i_{2}^{\prime
}+j_{2}^{\prime }=i_{2}+j_{2}$. 
\begin{equation}
0=\sum q^{i_{2}^{\prime }j_{1}^{\prime }}\lambda _{\left( i_{1}^{\prime
},i_{2}^{\prime }\right) }\mu _{\left( j_{1}^{\prime },j_{2}^{\prime
}\right) }  \label{prod eq}
\end{equation}%
It is easy to see $\left( i_{1}^{\prime },i_{2}^{\prime }\right) \preceq
\left( i_{1},i_{2}\right) $ or $\left( j_{1}^{\prime },j_{2}^{\prime
}\right) \preceq \left( j_{1},j_{2}\right) $ for each summand in equation %
\ref{prod eq}. \ By minimality $\lambda _{\left( i_{1}^{\prime
},i_{2}^{\prime }\right) }\mu _{\left( j_{1}^{\prime },j_{2}^{\prime
}\right) }=0$ unless $\left( i_{1}^{\prime },i_{2}^{\prime }\right) =\left(
i_{1},i_{2}\right) $ and $\left( j_{1}^{\prime },j_{2}^{\prime }\right)
=\left( j_{1},j_{2}\right) $. \ Now equation \ref{prod eq} becomes $q^{i_{2}
j_{1}}\lambda _{\left( i_{1},i_{2}\right) }\mu _{\left( j_{1},j_{2}\right)
}=0$, which is a contradiction since $q \neq 0$, $\lambda _{\left(
i_{1},i_{2}\right) } \neq 0$, and $\mu _{\left( j_{1},j_{2}\right) }\neq 0$.
\end{proof}

\section{Reducible Polynomials \label{Irreducible Sec}}

\begin{definition}
Let $f\in \mathcal{O}_{q}\left( F^{2}\right) \backslash \left\{ 0\right\} $.

\begin{enumerate}
\item We say $f$ is \emph{reducible} if $f=gh$ for some nonconstant $g,h\in 
\mathcal{O}_{q}\left( F^{2}\right) $. \ 

\item We say $f$ is \emph{irreducible} if it is nonconstant and not
reducible.

\item We call $f$ a \emph{quadratic form} if it can be written as in in
equation \ref{QF Eq}.%
\begin{equation}
f\left( x,y\right) =ax^{2}+bxy+cy^{2}  \label{QF Eq}
\end{equation}

\item The \emph{quantum discriminant} of $f$ in equation \ref{QF Eq}\ is $%
\Delta _{q}\left( f\right) =b^{2}-4acq$.
\end{enumerate}
\end{definition}

We use the quantum discriminant to provide a reducibility test for quadratic
forms, which we state as Theorem \ref{Reducible Theorem}.

\begin{theorem}
\label{Reducible Theorem}A quadratic form $f$ is reducible if and only if $%
\Delta _{q}\left( f\right) =d^{2}$ for some $d\in F$.
\end{theorem}

\begin{proof}
Write $f$ as in equation \ref{QF Eq}.\ \ We prove the result when $a\neq 0$
and $c\neq 0$. \ First suppose $f$ is reducible. By Lemma \ref{homogeneous
Lemma} we may write $f$ as in equation \ref{secondf} for some $\lambda ,\mu
\in F$. 
\begin{equation}
f\left( x,y\right) =a\left( x+\lambda y\right) \left( x+\mu y\right)
\label{secondf}
\end{equation}%
A straightforward check shows $\Delta _{q}\left( f\right) =d^{2}$\ where $%
d=a\left( \mu -q\lambda \right) $. To finish the proof we note there is a
factorization of $f$ by setting $\lambda =\dfrac{b+d}{2aq}$ and $\mu =\dfrac{%
b-d}{2a}$ in equation \ref{secondf}.
\end{proof}

\begin{example}
In $\mathcal{O}_{q}\left( \mathbb{R}^{2}\right) $ there are two distinct
factorizations of $x^{2}-y^{2}$ into irreducible polynomials if $q\ $is
positive and different from 1. 
\begin{eqnarray*}
x^{2}-y^{2} &=&\left( x+\frac{\sqrt{q}}{q}y\right) \left( x-\sqrt{q}y\right)
\\
&=&\left( x-\frac{\sqrt{q}}{q}y\right) \left( x+\sqrt{q}y\right)
\end{eqnarray*}%
If $q<0$ there are two distinct factorizations of $x^{2}+y^{2}$ in $\mathcal{%
O}_{q}\left( \mathbb{R}^{2}\right) $. \ Simply replace $x^{2}-y^{2}$ with $%
x^{2}+y^{2}$ and $\sqrt{q}$ with $\sqrt{-q}$ in the expression above.
\end{example}

\section{Factorization in Noncommutative Rings \label{Noncommutative Rings
Section}}

Throughout this section $R$ denotes an arbitrary (noncommutative) ring with
1. \ We cover some basic concepts in this more general setting.

\begin{definition}
\label{Ring Def}Let $r\in R$. \ 

\begin{enumerate}
\item We call $r$ \emph{normal} if for all $s,t\in R$ there exists $%
s^{\prime },t^{\prime }\in R$ such that $sr=rs^{\prime }$ and $rt=t^{\prime
}r$. \ 

\item We call $r$ \emph{central} if for all $s\in R$ we have $rs=sr$.

\item We say $r$ \emph{divides} $s\in R$ and write $r|s$ if there exists $%
t\in R$ such that $s=rt$ or $s=tr$.

\item We call $r$ \emph{prime} if $r|st$ implies $r|s$ or $r|t$ for all $%
s,t\in R$.

\item A two-sided ideal $I$ of $R$ is called \emph{completely prime} if $R/I$
is a domain.
\end{enumerate}
\end{definition}

Obviously every central element is normal.

\begin{remark}
\label{Center Description}In $\mathcal{O}_{q}\left( F^{2}\right) $ the
variables $x,y$ are normal. \ Prime elements of $\mathcal{O}_{q}\left(
F^{2}\right) $ are irreducible. \ If $q$ is a primitive $n^{\text{th}}$ root
of unity then $f\in \mathcal{O}_{q}\left( F^{2}\right) $ is central if and
only if the exponent of $x$ and the exponent of $y$ are both multiples of $n$
for every term in the support of $f$. \ If $q$ is not a root of unity then
the only central polynomials are scalars. \ These facts are easy to prove.
\end{remark}

The proof of Lemma \ref{normal ideal lemma} is straightforward.

\begin{lemma}
Suppose $r\in R$ is normal and $I$ is the smallest two-sided ideal
containing $r$. \label{normal ideal lemma}

\begin{enumerate}
\item $I=\left\{ sr:s\in R\right\} =\left\{ rt:t\in R\right\} $

\item $I$ is completely prime if and only if $r$ is prime
\end{enumerate}
\end{lemma}

We conclude this section with a brief account of factorization in
noncommutative rings. \ Prime elements in noncommutative Noetherian rings
were first defined by A. W. Chatters in \cite{Chatters}. \ Part 5 of
Definition \ref{Ring Def} reduces to Chatters' definition if $r$ is normal
and $R$ is a Noetherian prime polynomial identity ring. \ This follows
immediately from Lemma \ref{normal ideal lemma} and Jategaonkar's principal
ideal theorem (see \cite{Jategaonkar}).

\section{Prime Polynomials \label{Prime Poly Sec}}

Given a polynomial $f\left( x,y\right) \in \mathcal{O}_{q}\left(
F^{2}\right) $ and scalars $\lambda ,\mu \in F$ we may form a new polynomial
by replacing the variables $x$ and $y$ by the scalar multiples of $\lambda x$
and $\mu y$, respectively. \ We denote this new polynomial by $f\left(
\lambda x,\mu y\right) $.

\begin{lemma}
Suppose $p\left( x,y\right) \in \mathcal{O}_{q}\left( F^{2}\right) $ is
prime and $x^{i}y^{j}$ is in the support of $p\left( x,y\right) $ for some $%
i,j\in \mathbb{N}$. \ Then $p\left( x,y\right) =q^{i}p\left( x/q,y\right)
=q^{j}p\left( x,y/q\right) $.\label{technical}
\end{lemma}

\begin{proof}
The result is obvious if $p\left( x,y\right) $ is a scalar multiply of $x$
or $y$ so assume otherwise. \ We show $p\left( x,y\right) =q^{i}p\left(
x/q,y\right) $. \ The remaining equality can be proved similarly.

A straightforward check verifies the equation $yp\left( x/q,y\right)
=p\left( x,y\right) y$. \ This implies $p\left( x/q,y\right) $ is a scalar
multiple of $p\left( x,y\right) $ since $p$ is prime. \ Thus for some
nonzero scalar $\lambda $ we have $p\left( x,y\right) =\lambda p\left(
x/q,y\right) $. \ By comparing coefficients we obtain $\lambda =q^{i}$.
\end{proof}

Lemma \ref{technical} limits the number of possible prime polynomials in $%
\mathcal{O}_{q}\left( F^{2}\right) $. \ We explain how in the proof of
Theorem \ref{Main Theorem}.

\begin{theorem}
The elements $x,y$ of $\mathcal{O}_{q}\left( F^{2}\right) $ are prime. \label%
{Main Theorem}

\begin{enumerate}
\item If $q$ is not a root of unity then a prime polynomial must be a scalar
multiple of $x$ or a scalar multiple of of $y$.

\item Suppose $q$ is a primitive $n^{\text{th}}$ root of unity. \ A prime
polynomial is either central and irreducible or a scalar multiple of $x$ or
of $y$.
\end{enumerate}
\end{theorem}

\begin{proof}
We noted $x$ and $y$ are normal in remark \ref{Center Description}. Using
part 1 of Lemma \ref{normal ideal lemma} we may find ring isomorphisms $%
\mathcal{O}_{q}\left( F^{2}\right) /\left\langle x\right\rangle \cong F\left[
y\right] $ and $\mathcal{O}_{q}\left( F^{2}\right) /\left\langle
y\right\rangle \cong F\left[ x\right] $. Thus $x$ and $y$ are prime by part
2 of Lemma \ref{normal ideal lemma}.

Choose prime $p\left( x,y\right) \in \mathcal{O}_{q}\left( F^2 \right) $.
Then $p$ is irreducible by Remark \ref{Center Description}.

\begin{enumerate}
\item If $q$ is not a root of unity then Lemma \ref{technical} implies $%
p\left( x,y\right) =\mu x^{i}y^{j}$ for some $\mu \in F$ and $i,j \in 
\mathbb{N}$. Since $p$ is irreducible we conclude $i=0$ and $j=1$ or $i=1$
and $j=0$.

\item Suppose $q$ is a primitive $n^{\text{th}}$ root of unity and $p$ is
not a scalar multiple of $x$ or of $y$. Since $p$ is irreducible it cannot
be expressed as the product of a nonconstant polynomial with $x$ or with $y$%
. There must be terms $x^{i}y^{0}$ and $x^{0}y^{j}$ in the support $p$ for
some $i,j \in \mathbb{N}$. By applying Lemma \ref{technical} and comparing
coefficients we see the exponents of $x$ and $y$\ in every term in the
support of $p$ are multiples of $n$. Thus $p$ fits the description of
central elements in Remark \ref{Center Description}.
\end{enumerate}
\end{proof}

Theorems \ref{Test Theorem} and \ref{QF Prime} show prime polynomials exist
which are not scalar multiples of $x$ or $y$ even when $q\neq 1$.

\begin{theorem}
\label{Test Theorem}Suppose $p$ is an irreducible polynomial belonging to
either $F\left[ x\right] $ or $F\left[ y\right] $. \ If $p$ is central in $%
\mathcal{O}_{q}\left( F^{2}\right) $ then $p$ is prime in $\mathcal{O}%
_{q}\left( F^{2}\right) $.
\end{theorem}

\begin{proof}
The proof in both cases is similar, so we only prove the result when $p\in F%
\left[ x\right] $ is an irreducible and central polynomial in $\mathcal{O}%
_{q}\left( F^{2}\right) $. \ By remark \ref{Center Description} there is a
positive integer $n$ such that $q$ is a primitive $n^{\text{th}}$ root of
unity. \ We use results on skew polynomial rings, which can be found in \cite%
[Chapter 1]{GB} and \cite[Chapter 1]{GW}.

Set $R=F\left[ x\right] /\left( I\cap F\left[ x\right] \right) $. \ Then $R$
is a commutative domain, which we may view as an $F$-subalgebra of $\mathcal{%
O}_{q}\left( F^{2}\right) /I$. \ There is an $F$-algebra automorphism $%
\alpha :F\left[ x\right] \rightarrow F\left[ x\right] $ such that $\alpha
\left( x\right) =qx$ and $\alpha ^{n}=id_{R}$. \ Since $\alpha \left(
p\right) =p$ there is an induced $F$-algebra endomorphism from $R$ to $R$
which we also denote by $\alpha $. \ It is easy to show that $\alpha $ is an
automorphism of $R$. \ We can construct a skew-polynomial ring $R\left[ \hat{%
y};\alpha \right] $ and\ we will prove $\mathcal{O}_{q}\left( F^{2}\right)
/I\cong R\left[ \hat{y};\alpha \right] $.

By \cite[Excercise 1.1.I]{GW} there is an $F$-algebra homomorphism $\varphi :%
\mathcal{O}_{q}\left( F^{2}\right) \rightarrow R\left[ \hat{y};\alpha \right]
$ such that $\varphi \left( x\right) =x+I\cap R$ and $\varphi \left(
y\right) =\hat{y}$.\ \ Since $\varphi \left( p\right) =0$ there is an
induced $F$-algebra homomorphism from $\mathcal{O}_{q}\left( F^{2}\right) /I$
to $R\left[ \hat{y};\alpha \right] $, which we also denote by $\varphi $. On
the other hand, there is an $F$-algebra homomorphism $\psi :R\left[ \hat{y}%
;\alpha \right] \rightarrow \mathcal{O}_{q}\left( F^{2}\right) /I$ such that 
$\psi |_{R}=id_{R}$ and $\psi \left( \hat{y}\right) =y+I$ by \cite[Lemma
1.1.11]{GW}. \ A straightforward check shows $\psi $ and $\varphi $ are
inverses, thus $\mathcal{O}_{q}\left( F^{2}\right) /I$ is isomorphic to $R%
\left[ \hat{y};\alpha \right] $, which is a domain by \cite[Lemma 1.1.12]{GB}%
. \ Therefore $p$ is prime by Lemma \ref{normal ideal lemma}.
\end{proof}

\begin{theorem}
\label{QF Prime}Let $F$ be contained in $\mathbb{R}$ and consider a
quadratic form $f=ax^{2}+bxy+cy^{2}$ with $\Delta _{q}\left( f\right) <0$. \
Then $f$ is irreducible, but it may not be prime.

\begin{enumerate}
\item If $q\neq \pm 1$ then $f$ is not prime.

\item If $q=1$ then $f$ is prime.

\item If $q=-1$ then $f$ is prime if and only if $b=0$.
\end{enumerate}
\end{theorem}

\begin{proof}
First of all, $f$ is irreducible by Theorem \ref{Reducible Theorem}

\begin{enumerate}
\item If $q\neq \pm 1$ then $f$ is not prime by Theorem \ref{Main Theorem}.

\item If $q=1$ then $f$ is prime by \cite[Lemma 9.2.1]{BB}.

\item Suppose $q=-1$. \ If $b\neq 0$ then Remark \ref{Center Description}
implies $f$ is not central and hence not prime by Theorem \ref{Main Theorem}%
. We assume $b=0$ and show $f$ is prime. \ It is enough to handle the case
when $F=\mathbb{R}$ since $\mathcal{O}_{q}\left( F^{2}\right) $ is contained
in $\mathcal{O}_{q}\left( \mathbb{R}^{2}\right) $. Then $f=a\left(
x^{2}-\left( \lambda y\right) ^{2}\right) $ where $\lambda =\sqrt{-a^{-1}c}$.

Set $S=\mathcal{O}_{-1}\left( \mathbb{R}^{2}\right) /I$ where $I$ is the
principal ideal generated by $f$. \ Let $\mathbb{H}$ denote the division
ring of quaternions containing $i,j,k$ such that $i^{2}=j^{2}=k^{2}=-1$ and $%
ji=-ij$, $ki=-ik$, and $kj=-jk$. There is a ring homomorphism $\phi
:S\rightarrow \mathbb{H}\left[ t\right] $ determined by $\phi \left(
x+I\right) =it$ and $\phi \left( y+I\right) =\lambda ^{-1}jt$. \ 

We show that if $p\in \mathcal{O}_{-1}\left( \mathbb{R}^{2}\right) $ and $%
\phi \left( p+I\right) =0$ then $p\in I$. We may find $p_{0},p_{1}\in 
\mathbb{R}\left[ y\right] $ such that $p+I=\left( p_{0}\left( y\right)
+p_{1}\left( y\right) x\right) +I$ by applying the left division algorithm
in \cite[Theorem 2.1]{Effective} to \ $p$ divided by $f$. Equation \ref%
{phequation} is immediate. 
\begin{equation}
\phi \left( p+I\right) =\left( p_{0}\left( \lambda ^{-1}jt\right) \right)
+\left( tp_{1}\left( \lambda ^{-1}jt\right) \right) i  \label{phequation}
\end{equation}%
If $\phi \left( p+I\right) =0$ then we use equation \ref{phequation} to
obtain $p_{0}=p_{1}=0$ and $p\in I$. Thus $\phi $ is injective. Now we can
conclude that $S$ has no zero divisors since it is isomorphic to a subring
of $\mathbb{H}\left[ t\right] $. Thus $f$ is prime by Lemma \ref{normal
ideal lemma}.
\end{enumerate}
\end{proof}

\begin{example}
Theorems \ref{Test Theorem} and \ref{QF Prime} imply $x^{4}+2$ and $%
x^{2}-y^{2}$ are both prime in $\mathcal{O}_{-1}\left( \mathbb{Q}^{2}\right) 
$.
\end{example}

Now it is easy to see that $\mathcal{O}_{-1}\left( \mathbb{Q}^{2}\right) $
contains many prime polynomials, all of which are irreducible and central.
We finish with an example of a polynomial in $\mathcal{O}_{-1}\left( \mathbb{%
Q}^{2}\right) $ which is not prime even though it is irreducible and central.

\begin{example}
The polynomial $p=x^{4}+y^{4}\in \mathcal{O}_{-1}\left( \mathbb{Q}%
^{2}\right) $ is central. But $p$ is not prime since if $%
f=x^{3}-x^{2}y+xy^{2}-y^{3}$ and $g=x(1-x)-(1+x)y$ then $p$ divides $%
fg=(1-x)p$ but $p$ does not divide $f$ or $g$. \ We claim $p$ is also
irreducible. \ By Lemma \ref{homogeneous Lemma} it is enough to rule out the
following two cases. \ 

\begin{description}
\item[Case 1] $p=\left( a_{0}x^{3}+a_{1}x^{2}y+a_{2}xy^{2}+a_{3}y^{3}\right)
\left( b_{0}x+b_{1}y\right) $

\item[Case 2] $p=\left( a_{0}x^{2}+a_{1}xy+a_{2}y^{2}\right) \left(
b_{0}x^{2}+b_{1}xy+b_{2}y^{2}\right) $
\end{description}

The proof is similar in each case but the second is harder. \ We only
provide the details in case 2. \ Expanding and comparing coefficients yields
equations i-v, which we prove are inconsistent if $%
a_{0},a_{1},a_{2},b_{0},b_{1},b_{2}\in \mathbb{Q}$.%
\begin{equation*}
\begin{tabular}{rllrl}
\textrm{i} & $a_{0}b_{0}=1$ & \hspace{0.75in} & \textrm{iv} & $%
a_{1}b_{2}+a_{2}b_{1}=0$ \\ 
\textrm{ii} & $a_{0}b_{1}+a_{1}b_{0}=0$ &  & \textrm{v} & $a_{2}b_{2}=1$ \\ 
\textrm{iii} & $a_{0}b_{2}-a_{1}b_{1}+a_{2}b_{0}=0$ &  &  & 
\end{tabular}%
\end{equation*}%
If $a_{1}=0$ or $b_{1}=0$ then multiply iii by $b_{0}b_{2}$ and substitute i
and v. This gives $\left( b_{0}\right) ^{2}+\left( b_{2}\right) ^{2}=0$
which, together with i and v, yields a contradiction.

We are left with $a_{i}\neq 0$ and $b_{i}\neq 0$ for $0\leq i\leq 2$. \
Multiply ii by $b_{0}$ and substitute i to obtain equation \ref{b1 formula}. 
\begin{equation}
b_{1}+a_{1}\left( b_{0}\right) ^{2}=0  \label{b1 formula}
\end{equation}%
A similar calculation using iv and v gives $b_{1}+a_{1}\left( b_{2}\right)
^{2}=0$. \ Thus $\left( b_{0}\right) ^{2}=\left( b_{2}\right) ^{2}$ and $%
b_{0}=sb_{2}$ for $s=\pm 1$. \ Substitute $b_{0}=sb_{2}$ and $b_{2}=sb_{0}$
into iii and use i and v to obtain $2s=a_{1}b_{1}$. \ If we multiply
equation \ref{b1 formula} by $a_{1}$ and substitute $2s=a_{1}b_{1}$ we find $%
\left( a_{1}b_{0}\right) ^{2}=-2s$ which is a contradiction.
\end{example}

\noindent \textbf{Acknowledgements}

The authors wish to thank the referee for helpful comments and for
suggesting a simpler proof of Theorem \ref{Test Theorem}.

\end{document}